\documentclass[joc]{ipart-arxiv}

\RequirePackage[OT1]{fontenc}
\RequirePackage{amsthm,amsmath}
\RequirePackage[numbers,square]{natbib}
\RequirePackage[colorlinks]{hyperref}

\usepackage[english]{babel}
\usepackage[utf8]{inputenc}
\usepackage{amsfonts}
\usepackage{amsthm}
\usepackage{amssymb}
\usepackage{eqnarray,amsmath}
\usepackage{bbm}
\usepackage{bm}
\usepackage{url}

\startlocaldefs

\addto{\captionsenglish}{\renewcommand{\abstractname}{}}

\newcommand{\ac}{\text{ac}}

\theoremstyle{plain}

\newtheorem*{theorem-non}{Theorem}

\newtheorem{lemma}{Lemma}

\newtheorem*{prop-non}{Proposition}

\theoremstyle{definition}

\newtheorem*{defi-non}{Definition}
\newtheorem*{ex}{Example}

\theoremstyle{remark}

\endlocaldefs

\begin{document}

\begin{frontmatter}

\title{On sequences covering all rainbow $k$-progressions}
\runtitle{On sequences covering all rainbow $k$-progressions}


\begin{aug}
\author{\fnms{Leonardo} \snm{Alese}\thanksref{t1,t2}\ead[label=e1]{alese@tugraz.at}},
\address{Graz University of Technology\\ Institute of Geometry\\ Kopernikusgasse 24, 8010 Graz, Austria\\ \printead{e1}}
\author{\fnms{Stefan} \snm{Lendl}\thanksref{t1}\ead[label=e2]{lendl@math.tugraz.at}}
\address{Graz University of Technology\\ Institute of Discrete Mathematics\\ Steyrergasse 30, 8010 Graz, Austria\\ \printead{e2}}
\and
\author{\fnms{Paul} \snm{Tabatabai}\ead[label=e3]{tabatabai@protonmail.com}}
\address{Graz University of Technology\\ 8010 Graz\\ \printead{e3}}

\thankstext{t1}{The authors acknowledge the support of the Austrian Science Fund (FWF): W1230, Doctoral Program ``Discrete Mathematics''.}
\thankstext{t2}{The author acknowledges the support of SFB-Transregio 109
``Discretization in Geometry \& Dynamics'' funded by DFG and FWF (I 2978).}

\runauthor{Alese, Lendl, Tabatabai}
\end{aug}

    \begin{abstract}
        Let $\ac(n,k)$ denote the smallest positive integer 
        with the property that there exists an $n$-colouring $f$ of
        $\{1,\dots,\ac(n,k)\}$ such that for every $k$-subset
        $R \subseteq \{1, \dots, n\}$ there exists an (arithmetic)
        $k$\nobreakdash-progression $A$ in $\{1,\dots,\ac(n,k)\}$
        with $\{f(a) : a \in A\} = R$.

        Determining the behaviour of the function $\ac(n,k)$
        is a previously unstudied problem.
        We use the first moment method to give
        an asymptotic upper bound for $\ac(n,k)$ for the case $k = o(n^{1/{5}})$.
    \end{abstract}



\end{frontmatter}
\section{Introduction} 
\label{intro}

Let $a, k, d \in \mathbb{N}$. The set $A = \{a, a +d, a+2d, \dots, a+(k-1)d\}$
is called an (arithmetic) $k$-progression. We say $A$ has \emph{common difference} $d$.

Let $n, N \in \mathbb{N}$ ($n \leq N$) and let
$f:[N]\rightarrow [n]$ be an $n$-colouring of $[N]$. Let
$R \in \binom{[n]}{k}$ be a $k$-subset of $[n]$. We say a
$k$-progression $A$ in $[N]$ is \emph{$R$-coloured} if $\{ f(a) : a \in A \} = R$.
We call such a $k$-progression a \emph{rainbow $k$-progression}.
We say $f$ \emph{covers} $R$ if there is a $k$-progression in $[N]$ that is
$R$-coloured.
\begin{ex}
    The $6$-colouring
    $ f = (4,6,5,1,3,4,2,5,6,3,1,4)$
    of the interval $\{1,2,\dots,14\}$
    covers every $3$-subset of $\{1,\dots,6\}$; we give
    examples for some subsets:
    \begin{alignat*}{1}
        \{1,2,3\}\text{: } (4,6,5,\bm{1},3,4,\bm{2},5,6,\bm{3},1,4) \ \ \  \\[-0.35em]            
        \{3,4,5\}\text{: } (\bm{4},6,\bm{5},1,\bm{3},4,2,5,6,3,1,4) \ \ \  \\[-0.35em]            
        \{3,4,6\}\text{: } (\bm{4},6,5,1,\bm{3},4,2,5,\bm{6},3,1,4) \ \ \  \\[-0.35em]
        \{2,5,6\}\text{: } (4,6,5,1,3,4,\bm{2},\bm{5},\bm{6},3,1,4) \ \ \  \\[-0.35em]            
    \end{alignat*}
\end{ex}

For $n,k \in \mathbb{N}$ (where $k \leq n$), let $\ac(n,k)$ denote
the smallest positive integer such that there exists
an $n$-colouring $f$ of $[\ac(n,k)] = \{ 1,2,\dots, \ac(n,k) \}$ that covers
every $k$-subset of $[n]$.

Among related problems, the anti-van der Waerden numbers $\text{aw}([N],k)$ are well-studied in
Ramsey theory.
The number $\text{aw}([N],k)$ is defined to be the smallest positive integer $r$
such that every surjective $r$-colouring of $[N]$ contains at least one rainbow $k$-progression.

Butler~et~al.~\cite{butler} calculate
exact values of $\text{aw}([N],k)$ for small values of $N$ and $k$ and
give asymptotic results. 
Berikkyzy~et~al.~\cite{berikkyzy2016anti} give an exact formula for
$\text{aw}([N],3)$, proving a conjecture of Butler~et~al.~\cite{butler}.
Young~\cite{young2016rainbow} and Schulte~et~al.~\cite{young2018graph} study
generalizations of this problem to finite abelian groups and graphs, respectively.

The problem of studying anti-van der Waerden numbers
is about finding colourings avoiding all rainbow $k$-progressions.
Conversely, the problem we study in this work is about finding colourings
that \emph{do not} avoid \emph{any} rainbow $k$-progressions.

A wide range of problems about covering all $k$-subsets of $[n]$, on various
structures, are studied~\cite{chung1992universal, blackburn2012existence, cover}.

We prove the following asymptotic result.
\begin{theorem-non}
    \label{maintheorem}
    As $n$ tends to infinity, we have
    $$\ac(n,k) = \Omega\left(\sqrt{k\binom{n}{k}} \right).$$
    If $k = k(n) = o(n^{1/{5}})$, we have
    $$\ac(n,k) = \mathcal{O}\left(\log n \cdot
    e^{k/{2}}\cdot k^{-k/{2}+5/{4}}\cdot n^{k/{2}}\right)$$
\end{theorem-non}
Comparing the asymptotic upper and the asymptotic lower bound for the case 
$k~=~o(n^{1/{5}})$,
we see that the bounds differ by the factor $k \log n$.

The proof of the theorem is given in Section~\ref{secthmproof}.
The main tool of the proof (Lemma~\ref{mainlemma}) is shown in Section~\ref{lemmaproof}.
We achieve this by finding a lower bound on the expected number of
$k$-subsets of $[n]$ covered by a random colouring.

\section{Proof of Theorem}\label{secthmproof}

All asymptotics are to be understood with respect to $n$, where $n$ tends to infinity.

The lower bound in the theorem is a consequence of the fact that
an $n$-colouring of $[N]$ can only cover all $k$-subsets of $[n]$ if $[N]$
contains at least $\binom{n}{k}$ $k$-progressions.

The remainder of this section is dedicated to proving the upper bound given in the theorem.
To this end, as claimed let $k = k(n) = o(n^{1/{5}})$ and
$N = N(n) = \left \lceil  \sqrt{2}\sqrt{\frac{k-1}{k!}}\cdot n^{k/{2}} \right \rceil$.

The proof of the following lemma is given in Section~\ref{lemmaproof}.
\begin{lemma}
    \label{mainlemma}
        Let
        $\mathcal{F} \subseteq \binom{[n]}{k}$
        be a family of $k$-subsets of $[n]$. There exists an $n$-colouring $f^*$ of $[N]$
        such that the number of sets of $\mathcal{F}$ that are covered by $f^*$
        is at least
        $|\mathcal{F}|\left(\frac{1}{2} + o(1) \right).$
\end{lemma}

It follows that 
there exists an $n$-colouring $g_{0}$ of $[N]$
that covers
at least $\binom{n}{k}\left(\frac{1}{2} + o(1) \right)$
of the
sets of $\mathcal{F}_0 := \binom{[n]}{k}$.

Let $\mathcal{F}_1$ be the family of sets of $\mathcal{F}_0$
that have not been covered by $g_{0}$.
Applying Lemma~\ref{mainlemma} again, we obtain an $n$-colouring $g_{1}$ of $[N]$
that covers at least $|\mathcal{F}_1|\left(\frac{1}{2} + o(1) \right)$
of the sets of $\mathcal{F}_1$. We repeat this process $r$ times,
by defining $\mathcal{F}_i$ to be the family of $k$-subsets of $[n]$
not yet covered by any of the colourings $g_{0}, \dots,g_{i-1}$.

After $r$ iterations, the number of $k$-subsets of $[n]$
that are not covered by any of the constructed colourings is at most
$|\mathcal{F}_0| \left(\frac{1}{2} + o(1) \right)^r$. 
Setting $r = r(n,k) = \left \lceil \alpha \cdot k \log n \right \rceil$, where
$\alpha > \frac{1}{\log(2)}$, we get
$$|\mathcal{F}_0| \left(\frac{1}{2} + o(1) \right)^{r(n)} =
\binom{n}{k}\left(\frac{1}{2} + o(1) \right)^{r(n)} = o(1).$$
Thus, for sufficiently large $n$, after $r(n)$
iterations, every $k$-subset of $[n]$ is covered by at least one
of the colourings 
$$g_{0}, g_{1}, \dots, g_{r(n)-1}.$$

From the colourings $g_{0}, g_{1}, \dots, g_{r(n)-1}$
we construct an $n$-colouring $g$ of
$S := [r(n) \cdot N]$.
We split $S$ into $r(n)$
intervals of length $N$
and colour each of these intervals with the corresponding colouring $g_{i}$.
Formally, we set
$$g\left(i\cdot N+s\right) =
g_{i}(s) \ \ \ i \in \{ 0, \dots, r(n)-1 \}, \ s \in [N].$$
The colouring $g$ is an $n$-colouring of
$S~=~\left[\left \lceil \alpha\cdot k \log n \right \rceil
\cdot
\left \lceil \sqrt{2}\sqrt{\frac{k-1}{k!}}
\cdot n^{k/{2}} \right \rceil \right]$ 
that covers all $k$-subsets of $[n]$. It follows that
\begin{align*}
    \ac(n,k) = \mathcal{O}\left(k \cdot \log n \cdot
    \sqrt{\frac{k-1}{k!}} \cdot
    n^{k/{2}} \right).
\end{align*}
If $k = o(n^{1/{5}})$ tends to infinity as $n \rightarrow \infty$,
$$\ac(n,k) = \mathcal{O}\left(\log n \cdot
e^{k/{2}}\cdot k^{-k/{2}+5/{4}}\cdot n^{k/{2}}\right)$$
holds.

\section{Proof of Lemma \ref{mainlemma} using the probabilistic method}
\label{lemmaproof}

For $n, N, k \in \mathbb{N}$ (where $k \leq n \leq N$) let $f$ be a 
random $n$-colouring of $[N]$ (chosen uniformly at random from all such colourings).
For each $R \in \binom{[n]}{k}$
let $X_R$ be the indicator variable
of the event
\emph{``$f$ covers $R$''}.
Given a $k$-progression $A$ in $[N]$, let
$Y_{A,R}$ be the event
\emph{``The progression $A$ is $R$-coloured''}.

We are interested in the random variable
$\sum_{R \in \binom{[n]}{k}} X_R,$
which counts the number of $k$-subsets of $[n]$ that are covered by $f$.

For the sake of brevity, let $\text{AP}_k(N)$ denote the set of all $k$-progressions in $[N]$ and 
$\mathcal{H}_k(N) = \binom{\text{AP}_k(N)}{2}$ denote the set of all
unordered pairs of $k$\nobreakdash-progressions in $[N]$.
Note that $X_R$ is the indicator variable of 
the event $\bigcup\limits_{A \in \text{AP}_k(N)} Y_{A,R}$.

Using a Bonferroni inequality we obtain the following lower bound for $\mathbb{E}X_R$.
\begin{lemma}
\label{expression}
For every $k$-subset $R$ of $[n]$, the following holds:
\begin{align*}
&\mathbb{E}X_R = \mathbb{P}(X_R = 1) =
\mathbb{P}\bigg(\bigcup\limits_{A \in \text{AP}_k(N)} Y_{A,R}\bigg)
\\
& \geq \sum_{A \in \text{AP}_k(N)}\mathbb{P}(Y_{A,R})
- \sum_{\{ A,B \}  \in \mathcal{H}_k(N)}
\mathbb{P}\left(Y_{A,R} \cap Y_{B,R}\right) \\
& = \sum_{A \in \text{AP}_k(N)} \frac{k!}{n^{k}} - 
\sum_{i = 0}^{k-1} \sum_{\substack{\{ A,B \} \in \mathcal{H}_k(N) \\ |A \cap B| = i}}
\frac{k!(k-1)!}{n^{2k-i}}.
\end{align*}
\begin{flushright}
$\square$
\end{flushright}
\end{lemma}

To evaluate the lower bound from Lemma \ref{expression}, we need to count the
number $h(N,k) = |\text{AP}_k(N)|$ of $k$-progressions in $[N]$ and the numbers
$h_i(N,k)$, defined as the number of unordered pairs of $k$-progressions in
$[N]$ that intersect in exactly $i$ positions.\newpage

\begin{lemma}
\label{formulas}
As $N$ tends to infinity, the following asymptotic bounds hold:
\begin{itemize}
    \item $h(N,k) = \frac{N^2}{2k-2} + \mathcal{O}(N)$,
    \item $h_0(N,k) \leq \binom{h(N,k)}{2} = \frac{N^4}{8(k-1)^2} + \mathcal{O}(N^2/{k})$,
    \item $h_1(N,k) \leq h(N,k)k^2N = \mathcal{O}\left(N^3k\right)$,
    \item $h_j(N,k) \leq \binom{N}{2} \binom{\binom{k}{2}}{2} = \mathcal{O}\left(N^2k^4\right) \text{ for } j \geq 2.$
\end{itemize}
\begin{proof}
The formula for $h(N,k)$ is obtained by counting the number of ways to choose
the initial term and common difference of the progression.
We bound $h_0(N,k)$ by the number of unordered pairs of $k$-progressions.
The bound for $h_1(N,k)$ is obtained by fixing a $k$-progression and an element
of that progression; there are at most $kN$ $k$-progressions containing this element.
For each $j \geq 2$, $h_j(N,k)$ is bounded by the total number of pairs of $k$-progressions
intersecting in at least two positions.
For each pair of distinct elements there are at most $\binom{k}{2}$ $k$-progressions
containing both of them.

\end{proof}
\end{lemma}

We are ready to evaluate the lower bound from Lemma \ref{expression}.

\begin{lemma}
\label{longlemma}
    Let $k = k(n) = o(n^{1/{5}})$ and let 
    $N = N(n) = \left \lceil  \sqrt{2}\sqrt{\frac{k-1}{k!}}\cdot n^{k/{2}} \right \rceil$.
    Let $f$ be a random $n$-colouring of $[N]$. Then,
    for every $R \in \binom{[n]}{k}$ the inequality
    $$\mathbb{E}X_R \geq \frac{1}{2} + o(1)$$
    holds.
\begin{proof}
Using Lemma~\ref{expression} and 
the asymptotic bounds for $h$ and the $h_i$'s we get
\begin{align*}
\mathbb{E}X_R &\geq 
h(N) \frac{k!}{n^k} - h_0(N)\frac{k!k!}{n^{2k}} - h_1(N)\frac{k!(k-1)!}{n^{2k-1}}
- \sum_{i=2}^{k-1} h_i(N)\frac{k!(k-i)!}{n^{2k-i}} \\
 &\geq  \left(\frac{N^2}{2k-2} + \mathcal{O}(N) \right) \frac{k!}{n^k}
 - \left(\frac{N^4}{8(k-1)^2} + \mathcal{O}(N^2/{k}) \right)\frac{k!k!}{n^{2k}} \\
 &\qquad+  \mathcal{O}(N^3k)\frac{k!(k-1)!}{n^{2k-1}} +
 \mathcal{O}(N^2k^4) \sum_{i=2}^{k-1} \frac{k!(k-i)!}{n^{2k-i}} =: L(n).
\end{align*}
Only the terms $\frac{N^2}{2k-2}\frac{k!}{n^k}$
and $\frac{N^4}{8(k-1)^2}\frac{k!k!}{n^{2k}}$ are asymptotically relevant.
It follows from Stirling's formula that
$\mathcal{O}(N) \frac{k!}{n^k} = o(1)$,
$\mathcal{O}(N^2/{k})\frac{k!k!}{n^{2k}} = o(1)$, and
$\mathcal{O}(N^3k)\frac{k!(k-1)!}{n^{2k-1}} = o(1)$.
To see that
$\mathcal{O}(N^2k^4)\sum_{i=2}^{k-1}\frac{k!(k-i)!}{n^{2k-i}} = o(1)$,
we use the fact that 
the last term of the sum asymptotically dominates the sum of all other terms
and the assumption $k~=~o(n^{1/{5}})$

We are thus left with the following representation of $L(n)$:
$$ L(n) = \frac{N^2}{2k-2} \frac{k!}{n^k}
- \frac{N^4}{8(k-1)^2}\frac{k!k!}{n^{2k}} + o(1),$$
which, by our choice of $N$, gives $L(n) = \frac{1}{2} + o(1)$.

\end{proof}
\end{lemma}

Lemma~\ref{mainlemma} follows from Lemma~\ref{longlemma} by linearity of expectation.

\section{Conclusion}


Various generalizations of the problem we studied are possible, by replacing $[N]$
by another structure endowed with a sensible definition of $k$-progression. 
Structures of interest include cycles $\mathbb{Z}_{N}$, abelian groups and
graphs, which are already studied for anti-van der Waerden numbers.

\section*{Acknowledgements}

We would like to thank Christian Elsholtz for introducing us to topics and methods related to
this work. This paper grew out of research in the framework of the third
author's master's thesis.

\bibliography{lit}
\bibliographystyle{ieeetr}

\end{document}